\documentclass[11pt]{amsart}   

\usepackage{amscd,amssymb}
\newtheorem{theorem}{Theorem}[section]

\theoremstyle{definition}

\theoremstyle{remark}

\numberwithin{equation}{section}
\theoremstyle{corol}
\newtheorem{corol}[theorem]{Corollary} 

\begin{document}
\noindent                                             
\begin{picture}(150,36)                               
\put(5,20){\tiny{Submitted to}}                       
\put(5,7){\textbf{Topology Proceedings}}              
\put(0,0){\framebox(140,34){}}                        
\put(2,2){\framebox(136,30){}}                        
\end{picture}                                         

\vspace{0.5in}

\title[A short proof of a theorem of Morton Brown]
{A short proof of a theorem of Morton Brown on chains of cells}

\author{Vladimir Uspenskij}
\address{Department of Mathematics, 321 Morton Hall, Ohio
University, Athens, Ohio 45701, USA}

\email{uspensk@math.ohiou.edu}




\begin{abstract} 
Suppose that a topological space $X$ is 
the union of an increasing sequence of open subsets each of which 
is homeomorphic to the Euclidean space 
${\mathbf R}^n$. Then $X$ itself is homeomorphic to 
${\mathbf R}^n$. This is an old theorem of Morton Brown.
We observe that this theorem is an immediate consequence of other
two theorems of Morton Brown concerning near homeomorphisms and
cellular sets.
\end{abstract}

\maketitle

\def\a{\alpha}
\def\d{\delta}
\def\D{\Delta}
\def\g{\gamma}
\def\s{\sigma}
\def\Si{\Sigma}
\def\implies{\Rightarrow}
\def\R{{\mathbf R}}
\def\Rp{{\mathbf R_+^*}}
\def\o{\omega}
\def\O{\Omega}
\def\G{\Gamma}

\def\sB{{\mathcal B}}
\def\sC{{\mathcal C}}
\def\sE{{\mathcal E}}
\def\sF{{\mathcal F}}
\def\sG{{\mathcal G}}
\def\sH{{\mathcal H}}
\def\sJ{{\mathcal J}}
\def\sK{{\mathcal K}}
\def\sN{{\mathcal N}}
\def\sO{{\mathcal O}}
\def\sR{{\mathcal R}}
\def\sS{{\mathcal S}}
\def\sT{{\mathcal T}}
\def\sU{{\mathcal U}}
\def\sV{{\mathcal V}}

\def\sbs{\subset}
\def\rar{\rightarrow}
\def\e{\epsilon}

\def\ti{\times}
\def\obr{^{-1}}
\def\stm{\setminus}

\hyphenation{homeo-mor-phism}

\section{Introduction} 
Consider the following theorem due to Morton Brown \cite{B-cell}:

\begin{theorem}
Suppose that a topological space $X=\cup_{i=0}^\infty U_i$ is the union of
an increasing sequence of open subsets $U_i$ each of which is homeomorphic
to the Euclidean space $\R^n$. Then $X$ is homeomorphic to $\R^n$.
\label{main}
\end{theorem}

The aim of this paper is to give a very short proof of this theorem,
based on other two theorems by Morton Brown concerning near
homeomorphisms and cellular sets. These theorems read:

\begin{theorem}[\cite{B-near}, \cite{A}, {\cite[Theorem 6.7.4]{vM}}]
Let $(X_n)$ be an inverse sequence of compact metric spaces with limit
$X_\infty$. If all bonding maps $X_k\to X_n$ are near homeomorphisms, then so are
the limit projections $X_\infty\to X_n$.
\label{near}
\end{theorem}

\begin{theorem}%
[\cite{B-Sch}, {\cite[Theorem 5.2, Propositions 6.2 and 6.5]{D}}]
Let $F$ be a closed subset of the $n$-sphere $S^n$. The following conditions
are equivalent:
\begin{enumerate}
\item $F$ is cellular;
\item the quotient map $S^n\to S^n/F$ (which collapses $F$ to a point)
is a near homeomorphism;
\item the quotient space $S^n/F$ is homeomorphic to $S^n$.
\end{enumerate}
\label{cell}
\end{theorem}

\begin{corol}
Let $f:S^n\to S^n$ be a map of the $n$-sphere onto itself such that only one
point-inverse of $f$ has more than one point. Then $f$ is a near homeomorphism.
\label{cor}
\end{corol}

Let us explain the notions used in these theorems. A map $X\to Y$ between
compact spaces is a {\it near homeomorphism} if it is in the closure of the 
set of all homeomorphisms from $X$ onto $Y$, with respect to the compact-open
topology on the space $C(X,Y)$ of all maps from $X$ to $Y$. 
A (closed) 
{\it $n$-cell} is a space homeomorphic to the closed $n$-cube $[0,1]^n$.
A compact subset $C$
of a Hausdorff $n$-manifold $M$ is {\it cellular} 
if it has a base of open neighbourhoods in $M$ homeomorphic to $\R^n$, 
or, equivalently, if it 
is the intersection of a decreasing sequence $(B_k)$ of closed $n$-cells such that 
each $B_{k+1}$ lies in the interior of $B_k$.

Cellular sets were used in the beautiful paper \cite{B-Sch} to prove the Generalized
Schoenflies Theorem \cite[Theorem 6.6]{D}. For that, a stronger version of
Corollary~\ref{cor} was needed: every onto self-map of $S^n$ with two non-trivial
point-inverses is a near homeomorphism. This requires a little more effort.
For our purposes, the elementary Theorem~\ref{cell} suffices. To make the paper
less dependent on external sources, we show in Section~3 that Theorem~\ref{cell}
readily follows from Bing's Shrinking Criterion.

\section{A short proof of Theorem~\ref{main}}

The proof can be made one line: consider one point compactifications, and
apply Corollary~\ref{cor} and Theorem~\ref{near}. We now elaborate.

Let $X=\cup_{i=0}^\infty U_i$ be the union
of an increasing sequence of open subsets $U_i$ each of which is homeomorphic
to the Euclidean space $\R^n$. Note that $X$ must be Hausdorff: any two points
$x,y\in X$ lie in a Hausdorff open subspace $U_k$. Let 
$X_\infty=X\cup\{\infty\}$ be the one
point compactification of $X$. Let $F_i$ be the complement of $U_i$ in $X_\infty$.
Let $X_i=X_\infty/F_i$ be the space obtained by collapsing the closed set
$F_i$ to a point. Then $X_i$ is a one-point compactification of $U_i$ and hence
homeomorphic to the $n$-sphere $S^n$.

Since the sequence $(F_i)$ is decreasing, there are natural maps 
$p_i^j:X_j\to X_i$ for $j>i$, 
and we get an inverse sequence $(X_i)$ of $n$-spheres.
Since the quotient maps 
$p_i^\infty:X_\infty\to X_i$ separate points of $X_\infty$,
the limit of this sequence can be identified with $X_\infty$. 

The maps $p_i^j:X_j\to X_i$ have at most one non-trivial point-inverse. According
to Corollary~\ref{cor}, they are near homeomorphisms. In virtue of 
Theorem~\ref{near}, so is the map $p^\infty_0:X_\infty\to X_0$.
It follows that $X_\infty$ is homeomorphic to $S^n$. Hence $X$ is homeomorphic
to $\R^n$.

\section{Shrinkable decompositions and cellular sets}

To make the paper more self-contained, we show how to deduce Theorem~\ref{cell}
from Bing's Shrinking Criterion. 

A {\it decomposition} of a set is a cover by disjoint subsets.
If $G$ is a decomposition of $X$, a subset of $X$ is {\em $G$-saturated}
if it is the union of some elements of $G$.
A decomposition $G$ of a compact Hausdorff space $X$ is 
{\em upper semicontinuous} if one of the following equivalent conditions
holds: (1) there exists a compact Hausdorff space $Y$ and a continuous map
$f:X\to Y$ such that $G=\{f\obr(y):y\in Y\}$; 
(2) the set $\bigcup\{g\ti g: g\in G\}$ is closed in $X\ti X$;
(3) for every closed subset $F$ of $X$ its {\em $G$-saturation}
$\bigcup\{g\in G: g\hbox{ meets }F\}$ is closed. 
An upper semicontinuous decomposition
$G$ of a compact metric space $X$ is {\it shrinkable} 
if for every $\e>0$
and every cover $\sU$ of $X$ by $G$-saturated open sets 
there exists a homeomorphism $h$ of $X$ onto itself such that: (1) for every $g\in G$ the set
$h(g)$ has diameter $<\e$; (2) for every $x\in X$ there exists $U\in \sU$
such that $x\in U$ and $h(x)\in U$.

\smallskip

{\bf Bing's Shrinking Criterion 
(\cite[Theorem 5.2]{D}, \cite[Theorem 6.1.8]{vM}.} 
{\sl An onto map $f:X\to Y$ between compact
metric spaces is a near homeomorphism if and only if the decomposition
$\{f\obr(y):y\in Y\}$ of $X$ is shrinkable.}\qed

\smallskip

{\it Proof of Theorem~\ref{cell}}. 

$(1)\implies(2)$. If $F$ is a cellular set
in a compact $n$-manifold $M$, the decomposition $G_F$ of $M$ whose only 
non-singleton element is $F$ is shrinkable. This easily follows from the fact
that for every $\e>0$
there exists a homeomorphism of the $n$-cube $[0,1]^n$ onto itself which
is identity on the boundary and shrinks the subcube 
$[\e, 1-\e]^n$ to a set of small
diameter. Bing's Shrinking Criterion implies that the quotient map $M\to M/F$
is a near homeomorphism.

$(2)\implies (3)$ is trivial.

$(3)\implies (1)$. Suppose $S^n/F$ is homeomorphic to $S^n$. We want to prove
that $F$ is cellular. Let $U$ be an open neighbourhood of $F$. Denote the quotient
map $S^n\to S^n/F$ by $p$. Let $a\in S^n/F$ be the point onto which $F$ collapses,
$p(F)=\{a\}$. Then $p(U)$ is an open neighbourhood of $a$.
Since $S^n/F$ topologically is a sphere, there exists
a neighbourhood $V$ of $a$ such that $V\sbs p(U)$ and the complement $C$ of $V$
is $S^n/F$ is cellular. Then $p\obr(C)$ is cellular in $S^n$ 
(note that $p$ restricted to $S^n\stm F$ is a homeomorphism). 
From the first part of the proof (implication $(1)\implies(2)$) it follows that 
the complement of any cellular subset of $S^n$ is homeomorphic to $\R^n$. 
Thus $S^n\stm p\obr(C)=p\obr(V)$ is homeomorphic to $\R^n$,
and it is an open neighbourhood of $F$ which is contained 
in $U$. Since $U$ was arbitrary, $F$ is cellular.\qed

\bibliographystyle{amsplain}

\end{document}